\documentclass[12pt]{article}
\usepackage{amsfonts}
\usepackage{color}

\newtheorem{thm}{Theorem}[section]
\newtheorem{lemma}[thm]{Lemma}
\newtheorem{cor}[thm]{Corollary}

\newtheorem{prop}[thm]{Proposition}

\newtheorem{fact}[thm]{Fact}
\newtheorem{example}{Example}

\newcommand{\be}{\begin{equation}}
\newcommand{\ee}{\end{equation}}

\newcommand{\Fp}{\mathbb F _p}

\newcommand{\cF}{\mbox{$\mathcal F$}}

\newcommand{\cH}{\mbox{$\mathcal H$}}

\newcommand{\cM}{\mbox{$\mathcal M$}}

\newcommand{\cV}{\mbox{$\mathcal V$}}

\newcommand{\N}{\mathbb N}
\newcommand{\Z}{\mathbb Z}
\newcommand{\Q}{\mathbb Q}

\newcommand{\F}{\mathbb F}

\newcommand{\ve}[1]{\mathbf{#1}}

\newcommand{\monom}[2]{\ve{#1}^{\ve{#2}}}
\newcommand{\lex}{\mathrm{Lex}}

\newcommand{\stb}[2]{#1_1,\dots,#1_{#2}}

\begin{document}

\title{A note on linear Sperner families}
\author{G\'abor Heged\H {u}s\footnote{ \'Obuda University, Antal Bejczy Center for Intelligent Robotics, Kiscelli utca 82, Budapest, Hungary, H-1032,
hegedus.gabor@nik.uni-obuda.hu}
, Lajos R\'onyai\footnote{Institute of Computer Science ans Control, Hungarian Academy of Sciences; Department of Algebra, 
Budapest University of Technology, Budapest, lajos@info.ilab.sztaki.hu}
}

\footnotetext{Research supported in part by National Research, Development and
Innovation Office - NKFIH Grant No. K115288.}



\maketitle

\begin{center}
{\em This paper is dedicated to the memory of our teacher, 
colleague and friend, professor Tam\'as E. Schmidt.}
\end{center}

\begin{abstract}

In an earlier work we described Gr\"obner bases of the ideal of polynomials
over a field, which vanish on the set of characterisic vectors $\ve v\in \{0,1\}^n$
of the complete $d$ unifom set family over the ground set $[n]$.  In particular, it turns out that
the standard monomials of the above ideal are {\em ballot monomials}.
We give here a partial extension of this fact. A set family is a {\em linear Sperner system} 
if the characteristic vectors satisfy a linear equation 
$a_1v_1+\cdots +a_nv_n=k$, where the $a_i$ and $k$ are  positive integers. 
We prove that the lexicographic standard monomials for linear Sperner systems are also ballot monomials, provided that 
$0<a_1\leq a_2\leq \cdots \leq a_n$.  
As an application, we confirm a conjecture of Frankl in the special case of 
linear Sperner systems. 
\end{abstract}

\medskip


\medskip


\section{Introduction}

Throughout the paper $n$ will be a positive integer and $[n]$ stands 
for the set $\{1,2, \ldots, n\}$. The family of all subsets of $[n]$ is 
denoted by $2^{[n]}$.


Let $\F$ be a field.
$\F[x_1, \ldots, x_n]=\F[\ve x]$ denotes  the
ring of polynomials in commuting variables $x_1, \ldots, x_n$ over $\F$.
For a subset $F \subseteq [n]$ we write
$\ve x_F = \prod_{j \in F} x_j$.
In particular, $\ve x_{\emptyset}= 1$.

Let $\ve v_F\in \{0,1\}^n$ denote the
characteristic vector of a set
$F \subseteq [n]$.
For a family of subsets $\cF \subseteq 2^{[n]}$, let
$V(\cF) = \{\ve v_F : F \in \cF\} \subseteq \{0,1\}^n \subseteq \F^n$.
A polynomial $f\in \F[x_1,\ldots ,x_n]$ can be considered as a function
from $V(\cF)$ to $\F$ in the straightforward way.

Several interesting properties of finite set systems ${\mathcal
F}\subseteq 2^{[n]}$ can be formulated simply as statements
about  {\em polynomial functions on $V(\cF)$}. 
For instance, the rank of certain inclusion matrices 
can be studied in this setting 
(see for example Sections 2, 3 in \cite{MR}).
As for polynomial functions on $V(\cF)$,
it is natural to consider the ideal $I(V(\cF))$:
$$ I(V(\cF)):=\{f\in \F[\ve x]:~f(\ve v)=0 \mbox{ whenever } \ve v\in V(\cF)\}. $$
Substitution gives  an $\F$ algebra homomorphism from $\F[\ve x]$ to the $\F$ algebra  
of $\F$-valued functions on $V(\cF)$. 
A straightforward interpolation argument shows that this homomorphism is
surjective, and the kernel is exactly
$I(V(\cF))$. This way we can identify $\F[\ve x]/I(V(\cF))$ and the algebra of
$\F$ valued functions on $V(\cF)$. As a  consequence, we have 
\begin{equation} \label{function}
\dim _\F \F[\ve x]/I(V(\cF))=|\cF|.
\end{equation}


{\em Gr\"obner bases} and related structures of $I(V(\cF))$ were 
given for some families ${\cF}$, see \cite{MR} and the references therein. 
Before proceeding further, we recall some basic facts about to Gr\"obner 
bases and standard monomials. For details we refer to \cite{AL}, \cite{Bu}, \cite{CCS}, \cite{CLS}.

A linear order $\prec$ on the monomials over 
variables $x_1,x_2,\ldots, x_m$ is a {\em term order}, or {\em monomial
order}, if 1 is
the minimal element of $\prec$, and $\ve u \ve w\prec \ve v\ve w$ holds for any monomials
$\ve u,\ve v,\ve w$ with $\ve u\prec \ve v$. Two important term orders are the lexicographic
order $\prec_l$ and the deglex order $\prec _{d}$. We have
$$x_1^{i_1}x_2^{i_2}\cdots x_m^{i_m}\prec_l x_1^{j_1}x_2^{j_2}\cdots
x_m^{j_m}$$
iff $i_k<j_k$ holds for the smallest index $k$ such
that $i_k\not=j_k$. Concerning the deglex order, we have $\ve u\prec_{d} \ve v$ iff  either
$\deg \ve u <\deg \ve v$, or $\deg \ve  u =\deg \ve v$, and $\ve u\prec_l \ve v$.

The {\em leading monomial} ${\rm lm}(f)$
of a nonzero polynomial $f\in \F[\ve x]$ is the $\prec$-largest
monomial which appears with nonzero coefficient in the canonical form of $f$ as a linear
combination of monomials.

Let $I$ be an ideal of $\F[\ve x]$. A finite subset $G\subseteq I$ is a {\it
Gr\"obner basis} of $I$ if for every nonzero 
$f\in I$ there exists a $g\in G$ such
that ${\rm lm}(g)$ divides ${\rm lm}(f)$. In other words, the leading
monomials ${\rm lm}(g)$ for $g\in G $ generate the semigroup ideal of 
monomials
$\{ {\rm lm}(f):~f\in I\}$. It follows easily, that $G$ is
actually a basis of $I$, i.e. $G$ generates $I$ as an ideal of $\F[\ve x]$. 
A key fact is (cf. \cite[Chapter 1, Corollary
3.12]{CCS} or \cite[Corollary 1.6.5, Theorem 1.9.1]{AL}) that every
nonzero ideal $I$ of $\F[\ve x]$ has a Gr\"obner basis.

A monomial $\ve w\in \F[\ve x]$ is a {\it standard monomial for $I$} if
it is not a leading monomial for any $f\in I$. We denote by  ${\rm sm}(I)$ 
the set of standard monomials of $I$.
For a
nonzero ideal $I$ of $\F[\ve x]$ the set of monomials
${\rm sm}(I)$ is a downset: if $\ve w\in {\rm sm}(I)$, $\ve u,\ve v$ are
monomials from $\F[\ve x]$ such that $\ve w=\ve u \ve v$ then 
$\ve u\in {\rm sm}(I)$. Also,
${\rm sm}(I)$ gives a basis of the $\F$-vectorspace $\F[\ve x]/I$ in
the sense that  
every polynomial $g\in \F[\ve x]$ can be 
uniquely expressed as $h+f$ where $f\in I$ and
$h$ is a unique $\F$-linear combination of monomials from ${\rm sm}(I)$.

For a set family $\cF\subseteq 2^{[n]}$ the characteristic vectors 
in $V(\cF)$ are all 0,1-vectors, hence the polynomials $x_i^2-x_i$ 
all vanish  on $V(\cF)$. We infer that the standard monomials of 
$I(\cF):=I(V(\cF))$ are square-free monomials. Moreover, (\ref{function}) and the preceding 
paragraph imply that 

\begin{equation} \label{standard}
|\cF|=\dim_\F \F[\ve x]/I(\cF)= |{\rm sm}(I(\cF))|.
\end{equation}

Let $\ve a =(\stb an)$ be a vector with positive integer components $a_i$, and $k\in \N$.   
We define the set of vectors  $S(\ve a, k)\subseteq \{0,1\}^n\subseteq \F^n$ as follows:

$$ S(\ve a,k):=\{(v_1,\ldots ,v_n)\in{\{0,1\}}^n:~ \sum_{i=1}^n a_iv_i=k\}. $$

In this paper, with the exception of a brief remark, where
$\F=\F_p$ is considered, we
assume that $\F=\Q$. The set family corresponding to $S(\ve a,k)$ is a Sperner system or antichain. Sperner systems of the form 
$S(\ve a, k)$ are called {\em linear Sperner systems}. There are Sperner systems which are non linear. A simple example is the following family:

{\small

$$
T:=\{(1,1,0,0,0), (1,0,1,0,0), (1,0,0,1,0), (1,0,0,0,1), (0,1,1,0,0),(0,0,1,1,1)\}.
$$
}
Indeed, easy linear algebra shows that $S(\ve a,k)$ can contain the first 5 points of $T$ 
only if $a_1=a_2=\cdots =a_5$.

\medskip

The complete uniform family of all $d$ element subsets of $[n]$ is linear, in fact it is $S(\ve 1, d)$, where $\ve 1 =(1,\ldots, 1)$.
Following \cite{ARS}, in  \cite{HR1} we  described Gr\"obner bases and 
standard monomials for the ideals $I_{n,d}=I(S(\ve 1,d))$. 
Extensions and combinatorial applications 
were given in \cite{HR2}.


Assume that $\prec$ is an arbitrary term order  on $\F[\ve x]$ such that $x_1\succ x_2\succ \cdots \succ x_n.$
Let $0\leq d\leq n/2$ and denote by $\cM_{d,n}$ the set of all
monomials $\ve x_G$ such that $G=\{s_1<s_2<\ldots <s_j\}\subset [n]$
for which $j\leq d$ and $s_i\geq 2i$ holds for every $i$,
 $1\leq i\leq j$. These monomials $\ve x_G$ are the {\em ballot monomials} of degree at most $d$.
If $n$ is clear from the context, then we write $\cM_d$ instead of the more precise $\cM_{d,n}$. 
 It is known (see for example Lemma 2.3 and the following remark in \cite{ARS}) that 
$$ |\cM_d|={n \choose d}.$$
In \cite{ARS} it was also shown for the lex order $\prec_l$, and 
this was  extended in \cite{HR1}
to any term order $\prec$ such that $x_n\prec \cdots \prec x_1$, that
$\cM_d$ is the set of standard monomials for $I_{n,d}$ as well as
for $I_{n,n-d}$. Our main aim in this note is to prove a partial extension of the above result
to linear Sperner systems. Some of the results in \cite{BG} also served as motivation for our work in this direction.

\begin{thm} \label{main}
 Let $\ve a =(a_1,\ldots , a_n) \in  \Z^n$ be a vector such that $0<a_1\leq a_2\leq\ldots \leq a_n$, and 
$k$ be a natural number. Then the lexicographic standard monomials for $S(\ve a,k)$ are all ballot monomials. 
More precisely 
$$  {\rm sm}(I(S(\ve a,k))\subseteq \cM _{[n/2]}.$$    
\end{thm}

In the following example we give an explicit description of the lex standard monomials 
for $S(\ve a, k)$, when $a_1=\ldots =a_{n-1}=1$, and $a_n=t$ for some integer $t\geq1$.

\begin{example} \label{pelda}
Let $1\leq t\leq k\leq \frac {n-1}{2}$ be integers, $a_1=\ldots =a_{n-1}=1$, $a_n:=t$, and put $V:=S(\ve a, k)$. 
Then the set of the lex standard monomials of $I(V)$ is
$${\rm sm}(I(V))=\cM_{k,n-1}\cup \{  \ve mx_n:~ \ve m\in \cM_{k-t,n-1} \}.$$
\end{example}

The following fact is easy to see by symmetric chain decomposition (see Problem 13.20 in \cite{L}). Here we offer a
somewhat algebraic proof.

\begin{cor}
Suppose that the coordinates of $\ve a\in \Z^n$ are positive integers and $k\leq \frac n2$ is a natural number. Then 
$$
|S(\ve a,k)|\leq {n \choose k}. $$
\end{cor}

\noindent
{\bf Proof.} After possibly permuting the coordinates, we may assume that $0<a_1\leq a_2\leq  \cdots \leq a_n$. 
Observe also that a monomial $\ve x_G$ is a leading monomial for $I(S(\ve a,k))$ whenever  
$|G|>k$, hence 
$$ |S(\ve a,k)|=|{\rm sm}(I(S(\ve a,k))| \leq |\cM _{k}|={n \choose k}. $$
Here we first used (\ref{standard}), and the inequality follows from Theorem \ref{main}.  $\Box$

\medskip

A set family $\cF\subseteq 2^{[n]}$ {\em shatters} a subset $S\subseteq [n]$, if  for 
every $Y\subseteq S$ there exists an $F\in \cF$ such that $F\cap S=Y$. In \cite{F} Frankl conjectured that 
if a Sperner system  $\cF\subseteq 2^{[n]}$ does not shatter any $\ell$ element subset of $[n]$ for some integer 
$0\leq \ell \leq n/2$, then 
$$ |\cF|\leq {n \choose \ell-1}.$$
Here we confirm this conjecture for linear Sperner systems.

\begin{cor}
Suppose that the coordinates of $\ve a\in \Z^n$ are positive integers, $k,\ell$ are natural numbers, $\ell \leq n/2$, 
and $S(\ve a, k)$ does not shatter any $\ell$ element subset of $[n]$.  Then 
$$
|S(\ve a,k)|\leq {n \choose \ell-1}. $$
\end{cor}

\noindent
{\bf Proof.} After possibly permuting coordinates, we may again assume that  
$0<a_1\leq a_2\leq  \cdots \leq a_n$. By Theorem \ref{main} the lex standard monomials 
of $S(\ve a,k)$ are ballot monomials. Next we observe that the 
square-free monomials of degree at least $\ell$ are leading monomials for $S(\ve a,k)$. 
Indeed, let $S\subseteq [n]$ be a subset, $|S|\geq \ell$. Then $S$ 
is not shattered by $S(\ve a,k)$: there is a subset $Y\subseteq S$ such that no 
$F \subseteq [n]$ for which $\ve v_F \in S(\ve a,k)$ can give $Y=S\cap F$.  
Then the polynomial 
$$ f(\ve x)=\prod _{i\in Y}x_i\cdot \prod_{j\in S\setminus Y}(x_j-1) $$
vanishes on $S(\ve a,k)$ completely, and the leading monomial of $f$ is $\ve x_S$ (for an arbitrary term order). 
We obtain that  
$$  {\rm sm}(I(S(\ve a,k))\subseteq \cM _{\ell-1},$$
and hence 
     $$ |S(\ve a,k))|=     |{\rm sm}(I(S(\ve a,k))|\leq |\cM _{\ell-1}|= {n \choose \ell-1}.$$
$\Box$

\medskip

Let $p$ be a prime, $\ve a \in \N^n$ be a vector, $k\in \N$. We
consider the family 
$$ S_p(\ve a, k) =\{ \ve v \in{\{0,1\}}^n:~ \sum_{i=1}^n a_iv_i\equiv k (mod ~p) \}\subset \Q^n. $$ 
Note that $S_p(\ve a, k)$ is no longer a Sperner family. 
An interesting and useful fact is (see \cite{F2}, \cite{HR2})
that in degrees at most $p-1$ the deglex standard 
monomials for $S(\ve 1,k)$ and $S_p(\ve 1, k)$ are the same over 
$\F_p$. We have a similar but weaker statement for more 
general $\ve a$. Weaker in the sense that stronger upper bound 
is required for the degree of the monomials, and also in the sense that 
our argument works only for lex standard monomials
\footnote{A set $V\subseteq \{0,1\}^n$ can be considered as a subset of $\F^n$ 
for any field $\F$ . It is known that 
the set of lex standard monomials for $I(V)$ is independent 
of $\F$. This is  seen for example from 
Proposition \ref{gameprop}.}.

Let $t$ be an integer, $0<t\leq n/2$. We define $\cH _t$ as the set
of those subsets $\{s_1<s_2<\cdots <s_t\}$ of $[n]$ for which $t$
is the smallest index $j$ with  $s_j<2j$.

We have $\cH _1=\{\{1\}\}$, $\cH _2=\{\{2,3\}\}$, and $\cH
_3=\{\{2,4,5\},\{3,4,5\}\}$. It is clear that if $\{s_1< \ldots
<s_t\}\in \cH _t$, then $s_t=2t-1$, and $s_{t-1}=2t-2$ if $t>1$.

\begin{prop} \label{modp}
Suppose that $0< a_i\leq a_{i+1}$ for each $1\leq i\leq n-1$. Let $0<t\leq n/2$ be an integer, 
$T\in \cH_t$, and assume that $\sum_{i\in T} a_i<p$. Then  $\ve x_T$ is a lex leading monomial 
for  $S_p(\ve a, k)$. In particular, the conclusion holds when $\sum_{i\in [2t-1]} a_i\leq p$. 
\end{prop}

\medskip

In the next Section we prove Theorem \ref{main}, Proposition \ref{modp}, and discuss the details 
of Example \ref{pelda}.

\medskip


\section{Lex standard monomials for linear Sperner systems}

We shall need the following simple observations. 

\begin{fact} \label{teny} Let $G\subseteq [n]$. If the monomial $\ve x_G$ is not a 
ballot monomial, then there exists an integer $t>0$ and a $Y\in \cH_t$ such that 
$Y\subseteq G$. $\Box$ 
\end{fact}

\medskip

\begin{lemma}  \label{Ht} Let  $0<a_1\leq a_2\leq \cdots \leq a_n$ and $t$ be integers, $1<t\leq n/2$, $T\in \cH _t$.
Then 
$$ \sum_{i\in [2t-1]\setminus T}a_i\leq \sum_{i\in T\setminus \{2t-1\}}a_i < \sum_{i\in T}a_i. $$
\end{lemma}

\noindent
{\bf Proof.} We prove that there exists a bijective map $f$ from $T\setminus \{2t-1\}$ onto
$[2t-1]\setminus T$ such that $f(t)<t$ for every $t\in T\setminus \{2t-1\}$.

This holds because $T\in H_t$ and therefore if 
$$T\setminus \{2t-1\}=\{l_1<l_2<\cdots <l_{t-1}\},$$ 
then $l_i\geq 2i$ for $i=1,\ldots, t-1$. The map $f$ can be constructed 
inductively for $l_1,\ldots ,l_{t-1}$.

Indeed, we can set $f(l_1)=1$. Suppose now that we have constructed 
$f(l_j)$ for $j<i$. The numbers $l_j$ and $f(l_j)$ are all positive 
integers less than $l_i$ by the induction hypothesis. Their number is 
$2i-2$. In the interval $[1,2i-1]$ there are $2i-1$ integers, hence 
we have one, say $s$,  which is not  among the numbers considered 
previously. Then we can set 
$f(l_i)=s$.\footnote{An alternative way to construct $f$  is to 
observe first that if we write 
 $$\{1,2, \ldots, 2t-2\}=
\{l_1<l_2<\cdots <l_{t-1}\}\cup ^*\{s_1<s_2<\cdots
<s_{t-1}\}, $$
then we have $s_i<l_i$ for $i=1,\ldots, t-1$. We can then set $f(l_i)=s_i$ 
for every $i$.}

The existence of $f$ implies that 
$$ \sum_{i\in [2t-1]\setminus T}a_i\leq \sum_{i\in T\setminus \{2t-1\}}a_i   < \sum_{i\in T}a_i  . $$
This proves the lemma.
$\Box$

\medskip

Following \cite{FRR} and  \cite{MR}
we recall some facts about the Lex game, 
a method to determine
the lexicographic standard monomials of the vanishing ideal of a
finite set of points from $\F^n$, where $\F$ is an arbitrary field. 
Let $V\subseteq \F^n$ be a finite set, and $\ve w =
(\stb wn)\in \N^n$ an $n$ dimensional vector of natural numbers.
With these data as parameters, we define the Lex game $\lex(V;\ve
w)$, which is played by two players, Lea and Stan, as follows:

\noindent Both Lea and Stan know $V$ and $\ve w$. Their moves are:
\begin{enumerate}
\item[1 ] Lea chooses $w_n$ elements of $\F$. 
\item[]  Stan picks a value $y_n\in\F$, different from Lea's choices. 
\item[2 ] Lea now chooses $w_{n-1}$ elements of $\F$. 
\item[]  Stan picks a $y_{n-1}\in\F$, different from Lea's (last $w_{n-1}$) choices.
\item[\dots]  (The game proceeds in this way until the first
coordinate.) 
\item[$n$ ] Lea chooses $w_1$ elements of $\F$.
\item[]  Stan finally picks a $y_1\in\F$, different from Lea's
(last $w_1$) choices.
\end{enumerate}
The winner of the game is Stan, if in the course of the game he can select  a vector
$\ve y=(y_1,\dots,y_n)$ such that $\ve y\in V$, otherwise Lea wins
the game. If in any step there is no suitable choice $y_i$ for
Stan, then Lea wins also.

The game allows a characterization of the lexicographic
leading monomials and standard monomials for $V$ (Theorems 2 and 3 in \cite{FRR}).

\begin{prop}\label{gameprop}
Let $V\subseteq \F^n$ be a nonempty finite set and $\ve w\in \N^n$. Stan
wins $\lex(V;\ve w)$ if and only if $\monom xw $ is a lex standard monomial 
for $I(V)$. Equivalently, Lea wins the lex game if and only if  $\monom xw $ is a
lex leading monomial for the ideal $I(V)$.
\end{prop}

\medskip
\noindent
{\bf Proof of Theorem \ref{main}.} We may assume that $S=S(\ve a, k)$ is nonempty. 
By Fact \ref{teny} it suffices to prove that for any integer
$1\leq t\leq n/2$ and $T\in H_t$ the monomial $\ve x_T$ is a lexicographic 
leading monomial for $S$. Note that $|T|=t$ and $2t-1\in T$. 
The statement is clear for $t=1$, in fact $x_1$ is a leading monomial for 
$S$, because $a_1x_1+\cdots +a_nx_n-k$ vanishes on $S$. Suppose for the rest of the proof that $t>1$.

We employ the Lex game method, proving that Lea wins the 
the lex game $\lex(S,\ve v _T)$, where $\ve v_T$ is the 
characteristic vector of $T$. After Stan specifies the 
coordinate values $y_{2t},\ldots ,y_n$, what remains (if Lea has not won yet) is a lex game 
$Lex(V,\ve v_T)$ where $V\subseteq \{0,1\}^{2t-1}$ defined by 
$\sum _{i=1}^{2t-1}a_iv_i=k'$, for some positive integer $k'\leq k$, and  
$\ve v_T$ is viewed now as a vector in $\{0,1\}^{2t-1}$.

Let $\cV\subseteq 2^{[2t-1]}$  denote the set family whose corresponding set 
of characteristic vectors is $V$. 
We claim that $\cV$ does not shatter $T$. To be more specific, either there 
is no 
$F\in \cV$ such that $F\cap T=T$, or there is no
$G\in \cV$ such that $G\cap T=\emptyset$.
 
Suppose for contradiction that both $F, G \in \cV$ exist. Then 

\begin{equation} \label{egyes} 
\sum_{i\in [2t-1]\setminus T}a_i\geq \sum _{i\in G}a_i =k'
=\sum _{i\in F}a_i\geq \sum _{i\in T}a_i.
\end{equation}
But this is in contradiction with the inequality of Lemma \ref{Ht}, proving the claim. We obtained that $\ve x_T$ is 
a lex leading monomial for $V$, the corresponding vanishing polynomial being
either $\ve x_T$ or $\prod_{i\in T} (x_i-1)$. This implies, that Lea 
wins the game $\lex(V,\ve v_T)$, hence also $\lex(S,\ve v_T)$ as well. 
This finishes the proof. $\Box$

\medskip

\noindent
{\bf Remark.} We can exhibit a polynomial $Q(\ve x)\in \Q[\ve x]$ vanishing on $S$ with 
leading term $\ve x_T$ without using directly the Lex game method, as follows.
Let $U_0\subseteq \{0,1\}^{n-2t+1}$ denote the set of all vectors $(v_{2t},\ldots, v_n)$ which can 
be extended into a vector in $S$ which has 0 coordinate values 
everywhere in $T$.   
%
%
Let $P(x_{2t},\ldots ,x_n)\in \Q[\ve x]$ be a
polynomial which is 0 on $U_0$ and is 1 on $\{0,1\}^{n-2t+1}\setminus U_0$. 
Then set 
$$ Q(x_1,\ldots ,x_n)=\prod _{i\in T}(x_i-P(x_{2t},\ldots ,x_n)). $$
It is immediate that the lex leading term of $Q$ is $\ve x_T$, since $T\subseteq[2t-1]$.
Let $\mathbf v\in S$ be an arbitrary vector. On one hand, if 
$\mathbf u=(v_{2t},\ldots, v_n)\in U_0$, then 
$Q(\ve v)= \prod_{i\in T} v_i=0$ because by the claim in the preceding proof 
vectors from $U_0$ do not have extensions $\ve v\in S$ with $v_i=1$ for 
all $i\in T$. 
On the other hand, if $\ve u \in \{0,1\}^{n-2t+1}\setminus U_0$,
then $Q(\ve v)=\prod _{i\in T}(v_i-1)=0$ because $\ve u$ has 
no extension $\ve v\in S$ with values $v_i=0$ for all $i\in T$. 
We note also, that using the equality  $P^2=P$ of functions defined on $\{0,1\}^n$, we have 
$$Q(x_1,\ldots ,x_n)= \ve x_T+ \left(\prod _{i\in T}
(x_i-1)-\ve x_T\right)P(x_{2t},\ldots ,x_n), $$
again an equality of functions on $\{0,1\}^n$.  

\medskip
\noindent
{\bf Proof of Proposition \ref{modp}.}
The statement is clear for $t=1$. For a vector $\ve v\in S_p(\ve a, k)$ the value 
$v_1$ is determined by the rest of the values $v_i$  because $a_1$ is not 0 modulo $p$.
Henceforth we assume that $t>1$. As with Theorem \ref{main}, it suffices 
to show that a set 
$V\subseteq \{0,1\}^{2t-1}$ defined by 
$\sum _{i=1}^{2t-1}a_iv_i\equiv k'(mod ~p)$ for some integer $0\leq k'\leq p-1$,
can not shatter $T$. Assume the contrary. 
Let $\mathbf v= \mathbf v^{(0)} \in V$ be a vector  which is 0 at every coordinate from
$T$. Also let $\mathbf u=\mathbf v^{(t)} \in V$ be a vector which has
coordinates 1 at every coordinate from $T$. Using Lemma \ref{Ht}
we obtain  
$$ 0\leq \sum _{i\in [2t-1]} a_iv_i\leq \sum _{i\in [2t-1]\setminus T}a_i < \sum _{i\in T}a_i \leq \sum_{i\in [2t-1]} a_iu_i\leq \sum _{i\in [2t-1]}a_i<2p.$$
This is possible only if $\sum _i a_iv_i=k'$ and $\sum _i a_iu_i=k'+p$. 
Now for $\ell =1,\ldots ,t-1$ let $\mathbf v^{(\ell)} \in V$ be a vector which is 1 in
the first $\ell$ coordinates from $T$, and is 0 at the remaining $t-\ell$
coordinates belonging to $T$. It follows from the indirect hypothesis that such vectors $\ve v^{(\ell)} $ exist. 
The inequality $\sum _{i\in [2t-1]}a_i<2p$ implies that for every $\ell$ the sum 
$\sum _{i\in [2t-1]}a_iv^{(\ell)}_i$ is either $k'$ or $k'+p$.  
Clearly there must be an index $j$  with $0\leq
j<t$,  such that $\sum _{i\in [2t-1]} a_i v^{(j)}_i=k'$ and 
$\sum _{i\in [2t-1]} a_i v^{(j+1)}_i=k'+p$ . Set $\mathbf w =\mathbf
v^{(j+1)}-\mathbf v^{(j)}$. This vector has $\pm 1$ and 0 coordinates, moreover it 
is 0 on $T$ with the exception of $w_{s}=1$, where $s\in [2t-1]$ is the $(j+1)^{\rm th}$
element of $T$.  
Therefore we have  
\begin{equation} \label{prop07} 
p=\sum_{i=1}^{2t-1} a_iw_i\leq a_{s}+ \sum _{i\in [2t-1]\setminus T}a_i\leq a_s+\sum _{i\in T\setminus \{2t-1\}} a_i
\leq  \sum _{i\in T}a_i< p,
\end{equation}
a contradiction proving the statement. At the second inequality we used Lemma \ref{Ht} again, and 
$a_s\leq a_{2t-1}$ at the third. $\Box$ 

\medskip

\noindent
{\bf Verification of Example \ref{pelda}.}
We recall first the following recursion for the lex standard monomials (see  the proof of Theorem 4.3 in \cite{ARS}).
Let $V\subseteq \{0,1\}^n\subseteq {\F}^n$ be a subset of the Boolean cube. Define the sets of vectors 
$$
V_0:=\{\ve v\in \{0,1\}^{n-1}:~ (\ve v,0)\in V\}
$$
and
$$
V_1:=\{\ve v\in \{0,1\}^{n-1}:~ (\ve v,1)\in V\}.
$$
Then for the lex standard monomials of $I(V)$ we have 

$$
{\rm sm}(I(V))={\rm sm}(I(V_0))\cup {\rm sm}(I(V_1)) 
\cup \{ \ve m x_n:~ \ve m\in {\rm sm}(I(V_0))\cap {\rm sm}(I(V_1))\}.
$$
We apply this in the case  $V:=S(\ve a, k)$, $\ve a=(1,\ldots, 1, t)$. 
It is easy to see that

$$
V_0=\{(v_1,\ldots ,v_{n-1})\in \{0,1\}^{n-1}:~  \sum_{i=1}^{n-1} v_i=k\}
$$
and 
$$
V_1=\{(v_1,\ldots ,v_{n-1})\in \{0,1\}^{n-1}:~  \sum_{i=1}^{n-1} v_i=k-t\}.
$$
Observe that $V_0$ and $V_1$ correspond to complete uniform 
families. Then by the results of Section 2 and Theorem 4.3 of \cite{ARS} 
we have
$$
{\rm sm}(I(V_0))=\cM_{k,n-1}
$$
and
$$
{\rm sm}(I(V_1))=\cM_{k-t,n-1}.
$$
These together imply that
$$
{\rm sm}(I(V))=\cM_{k,n-1}\cup \cM_{k-t,n-1} \cup \{ \ve m x_n:~ \ve m\in \cM_{k,n-1}\cap \cM_{k-t,n-1} \}=
$$
$$
=\cM_{k,n-1}\cup \{ \ve m x_n:~ \ve m\in \cM_{k-t,n-1} \}.
$$
Here we used that $0\leq k-t<k\leq \frac {n-1}{2}$, and hence  $\cM_{k-t,n-1}\subseteq \cM_{k,n-1}$. $\Box$


\medskip



\medskip


\begin{thebibliography}{MM}


\bibitem{AL} Adams, W. W., Loustaunau, P.: An
Introduction to Gr\"obner bases. American Mathematical Society (1994)


\bibitem{ARS}Anstee, R.P., R\'onyai, L., Sali, A.: Shattering News.
Graphs and Combinatorics {\bf 18}, 59--73 (2002)

\bibitem{BF} Babai, L., Frankl, P.: Linear Algebra Methods in
Combinatorics with Applications to Geometry and Computer Science. The University of Chicago (1992)

\bibitem{BG} Balandraud \'E., Girard B.: A Nullstellensatz for Sequences Over $\Fp$.
Combinatorica {\bf 34}, 657--688 (2014)


\bibitem{Bu} Buchberger, B.: Gr\"obner-Bases: An Algorithmic Method in
Polynomial Ideal Theory. In: Bose, N.K. (ed.) Multidimensional
Systems Theory - Progress, Directions and Open Problems in
Multidimensional Systems Theory. 184-232,  Reidel Publishing
Company, Dodrecht - Boston - Lancaster (1985)


\bibitem{CCS}Cohen, A.M., Cuypers, H., Sterk, H. (eds.): Some
Tapas of Computer Algebra. Springer-Verlag, Berlin, Heidelberg (1999)

\bibitem{CLS} Cox, D., Little, J., O'Shea, D.:
Ideals, Varieties, and Algorithms. Springer-Verlag, Berlin,
Heidelberg (1992)

\bibitem{FRR} Felszeghy, B., R\'ath, B., R\'onyai, L.:
The lex game and some applications. J.~Symbolic Computation {\bf 41},
663--681 (2006)

\bibitem{F} Frankl, P.: Traces of antichains. 
Graphs Comb. {\bf  5}, 295--299 (1989)

\bibitem{F2} Frankl, P.: Intersection Theorems and mod $p$ Rank of Inclusion Matrices. Journal
of Combinatorial Theory, Series A. {\bf 54}, 85--94 (1990)

\bibitem{HR1} Heged\H us, G., R\'onyai, L.: Gr\"obner bases for complete
uniform families. J. of Algebraic Combinatorics {\bf 17}, 171--180 (2003)

\bibitem{HR2} Heged\H us, G., R\'onyai, L.:
Standard Monomials for $q$-uniform Families and a Conjecture of
Babai and Frankl.  Central European Journal of Mathematics. {\bf 1},
198--207 (2003)

\bibitem{L} Lov\'asz, L.: Combinatorial Problems and Exercises.
Akad\'emiai Kiad\'o, Budapest (1979)

\bibitem{MR} R\'onyai, L.,  M\'esz\'aros, T.:  Some combinatorial applications of Gr\"obner bases. 
In: Algebraic Informatics (Winkler, F. ed.), 4th International Conference, CAI
2011, Linz, Proceedings. 65--83,  Springer-Verlag, Heidelberg  (2011)






\end{thebibliography}
 \end{document}